\documentclass[12pt]{amsart}

\begin{document}

\title[$A$-hypergeometric functions]
{$A$-hypergeometric functions in transcendental questions\\ of algebraic geometry}
\author{A. V. Stoyanovsky}
\thanks{Partially supported by the grant RFBR 10-01-00536.}
\email{alexander.stoyanovsky@gmail.com}
\address{Russian State University of Humanities}

\begin{abstract}
We generalize the known constructions of $A$-hyp\-er\-geo\-metr\-ic functions.
In particular, we show that periods of middle dimension
on affine or projective complex algebraic varieties are $A$-hyp\-er\-geo\-metr\-ic functions of coefficients
of polynomial equations of these varieties.
\end{abstract}

\maketitle

{\bf 1. Introduction.}

In a series of papers, I. M. Gelfand with co-authors have introduced and studied the important class of
$A$-hyp\-er\-geo\-metr\-ic functions. The definition of $A$-hyp\-er\-geo\-metr\-ic
system of linear partial differential equations and a study of its solutions are given in [1].

The goal of this note is to provide a general algebro-geometric construction of $A$-hyp\-er\-geo\-metr\-ic
functions including the known constructions as particular cases. Let us recall these known constructions.

a) In [2] it is shown that periods of products of complex powers of arbitrary polynomials
of several variables, i.~e. integrals
\begin{equation}
\oint_C f_1(x_1,\ldots,x_m)^{\lambda_1}\ldots f_n(x_1,\ldots,x_m)^{\lambda_n}x_1^{\beta_1-1}\ldots x_m^{\beta_m-1}
dx_1\ldots dx_m
\end{equation}
over an $m$-dimensional real cycle $C$ with values in the corresponding local system,
are $A$-hypergeometric functions of the coefficients of the polynomials $f_1,\ldots,f_n$.

b) In [3] it is shown that
periods of exponent of an arbitrary polynomial, i.~e. integrals
\begin{equation}
\int_C e^{f(x_1,\ldots,x_m)}x_1^{\beta_1-1}\ldots x_m^{\beta_m-1}dx_1\ldots dx_m,
\end{equation}
where $C$ is a possibly non-compact $m$-dimensional contour (with values in the local system)
such that the expression under the integral
tends on it to zero at infinity, are $A$-hyp\-er\-geo\-metr\-ic functions
of the coefficients of the polynomial $f$.

c) Recall the fundamental theorem of B. Sturmfels [4] which is the ``constructive main theorem of algebra''.

{\bf Theorem 1.} {\it The complex roots of an arbitrary algebraic equation of one variable
\begin{equation}
f(x)=0
\end{equation}
form a
\emph(multivalued\emph) $A$-hypergeometric function of the coefficients of the polynomial $f(x)$.}

All these theorems are checked by direct differentiation showing that the required quantity (the integral or the root)
satisfies the $A$-hypergeometric system of partial differential
equations as a function of coefficients.

The present note arose in attempts to understand and to unify these constructions. The result
is a general construction
from the theory of periods of algebraic varieties, see the Main Theorem below. This theorem and its
corollaries, Theorems 2, 3 below, show that
the $A$-hyp\-er\-geo\-metr\-ic functions should play an important role in motivic constructions of algebraic geometry,
bridging the gap between purely analytic and purely algebraic theories.

The author thanks M. V. Finkelberg for clarifying discussions.

{\bf 2. Main theorem.}

{\bf Theorem.}
{\it Let
$$
x=(x_1,\ldots,x_m)\to z=f(x)=(f_1(x),\ldots,f_n(x))
$$
be a polynomial map of affine or projective spaces, where $f_i(x)=\sum_j a_{ij}x^j$, $1\le i\le n$,
$j=(j_1,\ldots,j_m)$, $x^j=x_1^{j_1}\ldots x_m^{j_m}$.
Let
$g=g(z_1,\ldots,z_n)$
be a \emph(multivalued in general\emph) holomorphic function,
which is quasi-homogeneous in the variables $z_i$ \emph(i.~e. $z_i\partial g/\partial z_i=\alpha_ig$
for certain number $\alpha_i$ and for each $i$\emph). Let
$\omega$ be a \emph(multivalued in general\emph) holomorphic differential $m$-form
of the variables $x$, which is quasi-homogeneous in the variables
$x_p$ \emph(i.~e. $x_p\partial\omega/\partial x_p=\beta_p\omega$\emph)\emph.
Then the integrals
\begin{equation}
\oint_C g(f(x))\omega
\end{equation}
are $A$-hyp\-er\-geo\-metr\-ic functions of the coefficients $a_{ij}$.
Here $C$ is an $m$-dimensional cycle with values in the local system determined by the multivalued
expression under the integral.
}

{\bf Proof.} The $A$-hyp\-er\-geo\-metr\-ic PDE's follow from the following computations:
\begin{equation}
\frac{\partial^q}{\partial a_{i_1j_{(1)}}\ldots\partial a_{i_qj_{(q)}}}\oint_C g(f(x))\omega=
\oint_C\frac{\partial^q g}{\partial z_{i_1}\ldots\partial z_{i_q}}(f(x))x^{j_{(1)}+\ldots+j_{(q)}}\omega;\\
\end{equation}
further, for each $p=1,\ldots,m$,
\begin{equation}
\begin{aligned}{}
&\sum_{i,j}j_pa_{ij}\frac{\partial}{\partial a_{ij}}\oint_C g(f(x))\omega
=\oint_C\sum_{i,j}\frac{\partial g}{\partial z_i}(f(x))
j_pa_{ij}x^j\omega\\
&=\oint_C x_p\frac{\partial g(f(x))}{\partial x_p}\omega=-\beta_p\oint_C g(f(x))\omega,
\end{aligned}
\end{equation}
since the form $\omega$ is quasi-homogeneous, and $g(f(x))\omega$ is closed; further, for each $i_0=1,\ldots,n$,
\begin{equation}
\begin{aligned}{}
&\sum_ja_{i_0j}\frac{\partial}{\partial a_{i_0j}}\oint_C g(f(x))\omega
=\oint_C\sum_j\frac{\partial g}{\partial z_{i_0}}(f(x))
a_{i_0j}x^j\omega\\
&=\oint_C \frac{\partial g}{\partial z_{i_0}}(f(x))z_{i_0}(x)\omega=\alpha_{i_0}\oint_C g(f(x))\omega,
\end{aligned}
\end{equation}
since the function $g(z)$ is quasi-homogeneous. $\square$

{\bf 3. Examples.} Example a) from the Introduction follows directly from the Main theorem if we put
$$
g(z_1,\ldots,z_n)=z_1^{\lambda_1}\ldots z_n^{\lambda_n},\ \
\omega=x_1^{\beta_1-1}\ldots x_m^{\beta_m-1}dx_1\wedge\ldots\wedge dx_m.
$$

Example b) follows by putting $n=1$,
$$
g(z)=e^z,\ \
\omega=x_1^{\beta_1-1}\ldots x_m^{\beta_m-1}dx_1\wedge\ldots\wedge dx_m.
$$
In this case,
equation (7) is not used (for equality (5) in this case implies more equations),
so that $g(z)$ does not need to be quasi-homogeneous.

Example c) follows by putting $m=n=1$, $g(z)=\frac1{2\pi i}\log z$, $\omega=dx$. Indeed, integrating by parts,
we have $\oint\log z\,dx=-\oint x dz/z$. Equation (7) is derived in the same way.

d) {\bf Theorem 2.} {\it Let
\begin{equation}
f(x_1,\ldots,x_k,y)=\sum c_{i_1\ldots i_kj}x_1^{i_1}\ldots x_k^{i_k}y^j=0
\end{equation}
be an arbitrary polynomial equation. Then
$y$ is a \emph(multivalued\emph) $A$-hypergeometric function of coefficients $c_{i_1\ldots i_kj}$,
depending on $x_1,\ldots,x_k$ as
on parameters.}

This theorem is obtained if we put $n=1$, $g(z)=\log z$, $\omega=dy$ in the Main Theorem.
In this case, the quasi-homogeneity
equation (6) is used only with respect to the variable $y$, and $x_1$,$\ldots$,$x_k$ play the role of parameters,
so that the form $g(f(x,y))\omega$ needs to be closed only with respect to $y$.

e) The following theorem is obtained if we put
$$
g(z_1,\ldots,z_n)=z_1^{\lambda_1}\ldots z_l^{\lambda_l}/(z_{l+1}\ldots z_n)
$$
in the Main Theorem and use the Cauchy residue theorem.

{\bf Theorem 3.} {\it The Gelfand--Leray integral
\begin{equation}
\begin{aligned}{}
\oint_C &f_1(x)^{\lambda_1}\ldots f_l(x)^{\lambda_l}x_1^{\beta_1-1}\ldots x_m^{\beta_m-1}\\
&dx_1\wedge\ldots\wedge dx_m/
(df_{l+1}\wedge\ldots\wedge df_n)|_{f_{l+1}(x)=\ldots=f_n(x)=0}
\end{aligned}
\end{equation}
over a cycle $C$ with values in the corresponding local system on the variety
\begin{equation}
f_{l+1}(x_1,\ldots,x_m)=\ldots=f_n(x_1,\ldots,x_m)=0,
\end{equation}
is an $A$-hypergeometric function of the coefficients of the polynomials $f_1$, $\ldots$, $f_n$.
}

In the case $m=n$, $l=0$, this Theorem is stated in [4].

{\it Remark.} All the results of this paper remain valid if one replaces affine or projective spaces
by arbitrary toric varieties.

\end{document}